\documentclass
[11pt,a4paper,oneside]{article}

\usepackage[english]{babel}

\usepackage{amsmath}
\usepackage{amsfonts}
\usepackage{amssymb}
\usepackage{amsthm}
\usepackage{makeidx}
\usepackage{fancyhdr}
\usepackage{showidx}
\usepackage{mathrsfs}
\usepackage{subfigure}
\usepackage{url}
\usepackage[dcucite]{harvard}

\providecommand{\abs}[1]{\left| #1 \right|}

\providecommand{\tonde}[1]{\left( #1 \right)}
\providecommand{\quadre}[1]{\left[ #1 \right]}
\providecommand{\graffe}[1]{\left\{ #1 \right\}}

\newcommand{\R}{\mathbb{R}}

\newcommand{\eps}{\varepsilon}

\newtheorem{teo}{Theorem}[section]

\theoremstyle{definition}

\theoremstyle{remark}
\newtheorem{oss}[teo]{Remark}

\title{Singular limit of a nonlinear fourth order inhomogeneous equation}
\author{C. Pocci}

\begin{document}
\maketitle
\begin{abstract}
In this paper we use the method of matched asymptotic expansions in order to obtain a geometric motion as the singular limit of a nonlinear fourth order inhomogeneous equation.
\end{abstract}
\textbf{Key words}: asymptotic expansions, reaction-diffusion equations, geometric evolution equations.\\
\emph{AMS Subject Classification}: 34E05, 35K57, 53C44.
\section{Introduction}
Geometric evolution equations are often studied as the singular limit of reaction-diffusion equations.
A typical example is the mean curvature flow equation 
\begin{equation}\label{AC_motion}
V=H,
\end{equation}
where $H$ is the curvature of the front and $V$ is the normal velocity to the front.
It is obtained as the singular limit of the Allen-Cahn equation \cite{chen}
\begin{equation}\label{AC_intro}
u_t=\Delta u-\frac{1}{\eps^2}\psi_u(u),\end{equation}
where $\eps>0$ is a small parameter and $\psi_u$ is the derivative of the double-well potential $\psi(u)=\frac12(1-u^2)^2$. If $\eps \to 0$, the solution $u$ tends to the minimal values of $\psi$, \emph{i. e.} $-1$ and $1$. Such a geometric evolution equation describes the motion of the boundary between regions where the limit equals $-1$ and $1$. In other words, the zero-level set of $u$ 
evolves according to \eqref{AC_motion} when the curvature flow has a smooth solution. Mean curvature flow appears in the description of interface evolutions in a variety of physical models; this is related to the property that such a flow is the gradient flow of the area functional and therefore occurs naturally in problems where a surface energy is minimized. As a matter of fact, equation \eqref{AC_intro} is the gradient flow, suitably rescaled, of the functional 
\begin{equation*}
\mathscr{F}_\eps(u) =
\int_\Omega\quadre{\frac{\eps}{2}\abs{\nabla
u}^2+\frac{1}{\eps}\psi(u)}dx,
\end{equation*}
where $\Omega$ is a bounded open subset of $\R^N$.
In \citeasnoun{LM}, the authors are interested in the following functional
\begin{equation}\label{LM_functional}
\mathscr{F}_\eps(u) = \int_\Omega\quadre{\frac{\eps}{2}\abs{\nabla
u}^2+ \frac{1}{\eps}\psi(u)}dx+\frac{1}{\eps}\int_\Omega\quadre{-\eps\Delta
u+\frac{1}{\eps}\psi_u(u)}^2 dx,
\end{equation}
which is related to a conjecture that De Giorgi made about the approximation, in the sense of $\Gamma$-convergence, of the Willmore functional, a functional depending on the curvature.
They study such a functional in the field of image segmentation, that is, the recovery of smooth boundaries in a picture.
The authors show that, in the case $N=3$, the zero level set of the gradient flow of \eqref{LM_functional}, suitably rescaled, approximates a front $\Gamma \subset \Omega$ that evolves according to the geometric motion
\begin{equation}\label{V_3D}
V=H-2\Delta_s H +4 H \tonde{K-\frac{H^2}{4}},
\end{equation}
where $V$ is the normal velocity of the interface, $\Delta_s H$ is the surface Laplacian, $H$ and $K$ are, respectively, the mean curvature and the Gaussian curvature of $\Gamma$ at $p(x)$, projection of $x$ on $\Gamma$ (\emph{i. e.} $H= k_1+k_2$ and $K=k_1\cdot k_2$, being $k_1$
and $k_2$ the principal curvatures of $\Gamma$). We observe that the term $\Delta_s H$ is related to the Cahn-Hilliard equation  \cite{cen}, since the zero level set of such an equation approximates, as $\eps \to 0$ and suitably rescaled, a front evolving with normal velocity proportional to the surface Laplacian of the mean curvature. \\
Due to the increasing interest in reaction-diffusion equations with spatially inhomogeneous reaction terms deriving from  \citeasnoun{nmhs} and  \citeasnoun{BL}, this paper is concerned with the asymptotic behaviour, as $\eps \to 0$, of the solution of the fourth order equation 
\begin{equation}\label{u_t_case_ab_app}
u_t=-\frac{1}{\eps}
w-\frac{2}{\eps^3}b^2(x)\psi_{uu}(u)w+\frac{4}{\eps} a(x)\nabla a(x)
\nabla w+\frac{2}{\eps} a^2(x)\Delta w,
\end{equation}
where
\begin{equation}\label{w_case_ab_app}
w= b^2(x)\psi_u(u)-2\eps^2 a(x) \nabla a(x) \nabla u-\eps^2
a^2(x)\Delta u.
\end{equation}
It is the gradient flow of the following functional
\begin{eqnarray}\label{func_ab}
\mathscr{F}_\eps(u) & =
&\int_\Omega\quadre{a^2(x)\frac{\eps}{2}\abs{\nabla
u}^2+b^2(x)\frac{1}{\eps}\psi(u)}dx\\
& + & \frac{1}{\eps}\int_\Omega\quadre{-\eps a^2(x)\Delta u-2\eps
a(x)\nabla a(x)\nabla u+b^2(x)\frac{1}{\eps}\psi_u(u)}^2
dx,\nonumber
\end{eqnarray}
where $a(x)$ and $b(x)$ are smooth and for some $a_0>0$, $a(x)\geq
a_0$ and $b(x)>0$  $\forall x \in \R^N$. In the same line of  \citeasnoun{LM} we show, by using the formal method of matched asymptotic expansions, that the zero level set of a solution of \eqref{u_t_case_ab_app}, \eqref{w_case_ab_app} approximates, as $\eps \to 0$ and rescaling the time $t \to \eps t$, a front evolving according to the geometric motion 
\begin{multline}\label{V_final}
V = a^2(x) H+4 H a^4(x)\tonde{K-\frac{H^2}{4}}\\+ 2 a^3(x) H^2 \nabla a(x) \nabla d+\frac{a^5(x)}{b(x)}H\Delta\tonde{\frac{b(x)}{a(x)}}-\frac{2}{15}\pi^2 a^6(x) H \quadre{\nabla \tonde{\frac{b(x)}{a(x)}}}^2\\
- a^2(x) A+2a^2(x) B+a(x)b(x)C\tonde{\frac{2}{15}\pi^2-1}-\frac{\pi^2}{15} a(x) b(x) D,
\end{multline}
in the case $N=3$, where $\nabla d$ and $A, B, C, D$ are specified in sections \ref{formulation}, \ref{calculations}. The paper is organized as follows. In section \ref{MAE}, we briefly describe the idea of the method of matched asymptotic expansions. Section \ref{formulation} is devoted to the introduction of the problem and the instruments that we use to solve it. Finally, Section \ref{calculations} contains the main calculations for the formal derivation of \eqref{V_final}.
\section{The method of matched asymptotic expansions}\label{MAE}
In this section, we introduce the idea of the formal tool that we will adopt throughout the paper. Usually, the method of matched asymptotic expansions is used to solve singular perturbations problems. 
The general method with perturbation problems is to seek an expansion with
respect to the asymptotic sequence $\graffe{1, \eps, \eps^2, \ldots}$, as $\eps \to 0$
\begin{equation*}
u^{\eps}\sim U_0+\eps U_1+\eps^2 U_2+\ldots, \quad \eps \to 0,
\end{equation*}
for functions $U_0, U_1, U_2, \ldots$ which have to be determined. In a singular perturbation problem, the regular methods produce an expansion that fails, at some point, to be valid over the whole domain. Then, the idea of the method is the following. The domain can be divided into two subdomains. On one subdomain, a solution is accurately approximated by an asymptotic series, obtained by treating the problem as a regular perturbation; we refer to it as the outer solution. On the other subdomain, the approximation cannot be accurate, since the perturbation terms in the problem are not negligible there. An approximation in the form of an asymptotic series is obtained there by handling this part of the domain
as a separate perturbation problem; we call it the inner solution. When the validity regions of the outer and inner
expansions overlap, the outer and inner solutions are combined through a certain process (matching) in such a way
that an approximate solution for the complete domain is found. 
\section{Formulation of the problem and preliminaries}\label{formulation}
The gradient flow of the functional \eqref{func_ab} is given by equations \eqref{u_t_case_ab_app}, \eqref{w_case_ab_app}.
We rescale the time $t_1=\eps t$.  The boundary conditions are obtained by taking the first variation of the functional (\ref{func_ab}) and provide that the first variation does not contain boundary terms. The problem to solve is then
\begin{equation}\label{problem_ab}
\left\{\begin{array}{ll} 

\vspace{1mm}

 \eps^4u_{t_1}=-\eps^2
w-2 b^2(x)\psi_{uu}(u)w+4\eps^2 a(x)\nabla a(x) \nabla w\\

\vspace{3mm}

+2\eps^2 a^2(x)\Delta w, \quad ${\rm in}$\quad\Omega \times (0, T)\\

\vspace{1mm}

w=b^2(x)\psi_u(u)-2\eps^2 a(x)\nabla a(x)\nabla u\\

\vspace{3mm}

-\eps^2 a^2(x) \Delta u, \quad ${\rm in}$\quad\Omega \times (0, T)\\

\vspace{3mm}

w=0, \quad ${\rm on}$\quad\partial\Omega\\

\vspace{3mm}

\overline n \cdot \quadre{\eps a^2(x) \nabla u+\frac{2}{\eps} a^2(x) \nabla
w}=0, \quad ${\rm on}$\quad\partial\Omega
\end{array} \right.,
\end{equation}
$\overline n$ being the outward unit normal to $\partial \Omega$.   
We denote the solutions of \eqref{u_t_case_ab_app} and \eqref{w_case_ab_app} with $u$ and $w$ respectively, keeping in mind their dependance on $\eps$. Moreover, we indicate the interface as
\begin{equation*}
\Gamma=\cup_{t\geq 0}\tonde{\Gamma_t \times \graffe{t}},
\end{equation*}
where $\Gamma_t=\graffe{x \in \R^N: u(x, t)=0}$. 
We assume that the zero level set of a solution $u$ of \eqref{u_t_case_ab_app} divides the domain $\Omega$ into two subdomains, $\Omega_+$ and $\Omega_-$.
Following \citeasnoun{cen}, \citeasnoun{pego} and \citeasnoun{LM}, we assume 
that $u$ and $w$ have the expansions
\begin{equation}\label{out_exp_u}
u(x, t_1)\sim u^0+\eps u^1+\eps^2 u^2+\ldots
\end{equation}
and 
\begin{equation}\label{out_exp_w}
w(x, t_1)\sim w^0+\eps w^1+\eps^2 w^2+\ldots
\end{equation}
away from the interface $\Gamma$. 
We refer to \eqref{out_exp_u} and \eqref{out_exp_w} as the outer expansions.
The construction of an inner solution for $x$ near $\Gamma$ is more complicated: indeed, a correct scale must be decided in order to get inner expansions which match well the outer expansions.
For this reason, in a small neighborhood of $\Gamma$, a stretched normal distance to the front is defined
\begin{equation}\label{z}
z=\frac{d(x, t_1)}{\eps},
\end{equation}
being $d(x, t_1)$ the signed distance from the point $x \in \Omega$ to $\Gamma_{t_1}$, such that $d>0$ for $x \in \Omega_+$ and $d<0$ for $x \in \Omega_-$.
An ansazt is introduced as follows for the inner expansions
\begin{equation}\label{inn_exp_u}
u(x, t_1)= U(z, x, t_1)\sim U^0+\eps U^1+\eps^2 U^2+\ldots,
\end{equation}
\begin{equation}\label{inn_exp_w}
w(x, t_1)= W(z, x, t_1)\sim W^0+\eps W^1+\eps^2 W^2+\ldots.
\end{equation}
As in \citeasnoun{pego} and \citeasnoun{LM}, it is required that the quantities depending on $(z, x, t_1)$ do not change when $x$ varies normal to $\Gamma$, keeping $z$ fixed.
Then, we have the following setting
\begin{equation*} 
V=\partial_{t_1} d(x, t_1),
\end{equation*}
being $V$ the normal velocity of $\Gamma$ in the $t_1$ timescale, positive when the front moves toward $\Omega_-$.
Furthermore, indicating with $\nabla d(x, t_1)$ the unit normal to $\Gamma$ pointing toward $\Omega_+$, we have
\begin{eqnarray*}
\nabla u & = & \frac{1}{\eps} \nabla d U_z+\nabla_x U \nonumber\\
& = & \nabla_x U^0+\eps \nabla_x U^1+\eps^2 \nabla_x U^2+\frac{1}{\eps}\nabla d(U^0_z+\eps U^1_z+\eps^2 U^2_z)+\ldots,
\end{eqnarray*}
\begin{eqnarray*}
\Delta u & = & \Delta_x U+\frac{1}{\eps}U_{z}\Delta
d+\frac{1}{\eps^2}U_{zz} \nonumber \\
& = & \frac{1}{\eps^2}U^0_{zz}+\frac{1}{\eps}(U^0_z \Delta d+U^1_{zz})\ldots,
\end{eqnarray*}
\begin{equation*}
\partial_{t_1}u=\partial_{t_1}U+\frac{1}{\eps}VU_z;
\end{equation*}
analogous equations are true for $w$ and $W$.
From \citeasnoun{gilbarg_trudinger}, in the case $N=3$
\begin{equation*}
\Delta d(x, t_1)=\frac{k_1(p(x))}{1+k_1(p(x)) d(x, t_1)}+\frac{k_2(p(x))}{1+k_2(p(x)) d(x, t_1)},
\end{equation*}
with $p(x)$ projection of $x$ on $\Gamma$.
However, from \eqref{z} and expanding on powers of $\eps$, it is easy to get
\begin{equation}\label{exp_d}
\Delta d(x, t_1) = H-\eps z(H^2-2K)+ O(\eps^2),
\end{equation} 
where $H$ and $K$ denote respectively the mean curvature and the gaussian curvature of the interface at $p(x)$.
From the definition of zero level set of $\Gamma$, the normalization conditions hold
\begin{equation}\label{normalization_ab}
\quad U^k(0, x, t_1)=0, \quad k=0, 1, \ldots.
\end{equation}
Finally, using the notation $u^i_{\pm}(x, t_1)=\lim_{r\to 0^{\pm}}u^i(x+r\nabla d(x, t_1), t_1)$, the following matching conditions are required:
\begin{equation}\label{match1}
u^0_{\pm}(x, t_1)=\lim_{z\to\pm \infty}U^0(z, x, t_1),
\end{equation}
\begin{equation}\label{match2}
\lim_{z\to \pm \infty}(u^1_{\pm}+z\nabla d\nabla u^0_{\pm})(x, t_1)=\lim_{z\to \pm \infty} U^1(z, x, t_1),
\end{equation}
\begin{equation}\label{match3}
\lim_{z\to \pm \infty}(u^2_{\pm}+z\nabla d \nabla u^1_{\pm}+\frac12 z^2 D^2 u^0_{\pm})(x, t_1)=\lim_{z\to \pm \infty} U^2(z, x, t_1),
\end{equation}
where $D$ denotes the directional derivative along $\nabla d(x, t_1)$ (see \citeasnoun{LM} for further details). The expansions of $w$ and $W$ must satisfy analogous matching conditions.
\section{Asymptotics}\label{calculations}
In this section, we present a formal derivation of the interface equation \eqref{V_final}.
\subsection{The outer solution}
At the beginning, we start with the solution far from $\Gamma$.
We notice that the minimizers of \eqref{func_ab} are expected to satisfy $u(x)=\pm 1$ almost everywhere; for this reason, we seek a solution $u$ of \eqref{u_t_case_ab_app} which is equal to $-1$ and $1$ in $\Omega_-$ and $\Omega_+$ respectively, in the timescale $t_1=\eps t$, as $\eps$ approaching zero.
The boundary conditions on $\partial \Omega$
yield
\begin{eqnarray*}
w^i& = & 0, \quad i=0, 1, \ldots,\\
\overline n \cdot a^2(x)\nabla w^0 = \overline n \cdot a^2(x)\nabla w^1 & = & 0,\\
\overline n \cdot(a^2(x)\nabla u^i+2a^2(x)\nabla w^{i+2}) & = & 0, \quad
i=0, 1, \ldots.
\end{eqnarray*}
The substitution of the outer expansions \eqref{out_exp_u} and \eqref{out_exp_w} into  (\ref{problem_ab}) leads at zero order in $\eps$ to the following problem
\begin{equation}\label{out_zero_case_ab}
\left\{\begin{array}{ll} b^2(x)\psi_{uu}(u^0)w^0=0 \\
w^0=b^2(x)\psi_u(u^0)\\
w^0=0, \quad {\rm on}\quad \partial \Omega\\
\overline n \cdot a^2(x) \nabla w^0=0, \quad {\rm on}\quad \partial \Omega
\end{array} \right..
\end{equation}
A solution of problem \eqref{out_zero_case_ab}, in accordance with our expectation, is 
\begin{equation*}
u^0(x, t_1)=\left\{\begin{array}{ll} +1,\quad x \in \Omega_+\\
-1, \quad x \in \Omega_-
\end{array} \right., \quad w^0(x, t_1)=0 \quad {\rm in }\quad \Omega.
\end{equation*}
Moreover, at first order in $\eps$, we have
\begin{equation}\label{out_one_case_ab}
\left\{\begin{array}{ll} b^2(x)\quadre{\psi_{uu}(u^0)w^1+\psi_{uuu}(u^0)u^1 w^0}=0 \\
w^1=b^2(x)\psi_{uu}(u^0)u^1\\
w^1=0, \quad {\rm on}\quad \partial \Omega\\
\overline n \cdot a^2(x) \nabla w^1=0, \quad {\rm on}\quad \partial \Omega
\end{array} \right.,
\end{equation}
whose solution is given by
\begin{equation*}
u^1(x, t_1)=0, \quad w^1(x, t_1)=0.
\end{equation*}
Going further in the analysis of second and third order, we have that $u^2$ and $u^3$ are null, hence the solution in the outer region is:
\begin{equation*}
u(x, t_1)=\pm 1+O(\eps^4).
\end{equation*}
\subsection{The inner solution}
In the following, we deal with the solution near $\Gamma$. We substitute expansions \eqref{inn_exp_u} and \eqref{inn_exp_w} into the equations of problem \eqref{problem_ab}.
\subsubsection*{Zero order}
As already done in \citeasnoun{cen} and \citeasnoun{LM}, we assume $W^0(z, x, t_1)=0$. At zero order in $\eps$ the equation (\ref{w_case_ab_app}) becomes:
\begin{equation}\label{inn_zero_case_ab}
a^2(x)U^0_{zz}-b^2(x)\psi_u(U^0)=0.
\end{equation}
The unique monotone increasing solution of (\ref{inn_zero_case_ab}) in accordance with the normalization condition (\ref{normalization_ab}) and the matching condition (\ref{match1}), is
\begin{equation}
U^0(z, x, t_1)=\tanh\quadre{\frac{b(x)}{a(x)}z}.
\end{equation}
\begin{oss}
We notice that the solution $U^0$ satisfies
\begin{equation}\label{relazione}
U^0_{zz}=-2 \frac{b(x)}{a(x)}U^0 U^0_z.
\end{equation}
\end{oss}
\subsubsection*{First order}
At first order in $\eps$, we have
\begin{equation}
a^2(x)W^1_{zz}-b^2(x)\psi_{uu}(U^0)W^1=0,
\end{equation}
\begin{equation}\label{inn_one_case_ab2}
W^1=b^2(x)\psi_{uu}(U^0)U^1-a^2(x) U^0_z H-a^2(x)U^1_{zz}.
\end{equation}
Let us define the operator $L$:
\begin{equation*}
LW^1=a^2(x)W^1_{zz}-b^2(x)\psi_{uu}(U^0)W^1.
\end{equation*}
Clearly, from (\ref{inn_zero_case_ab}), we obtain
\begin{equation*}
LU^0_z=0,
\end{equation*}
therefore we can assume
\begin{equation*}
W^1(z, x, t_1)=\mu(x, t_1)U^0_z\left[\frac{b(x)}{a(x)} z\right],
\end{equation*}
with $\mu(x, t_1)$ a bounded function to be determined. In this setting, equation (\ref{inn_one_case_ab2}) becomes
\begin{equation}\label{LU1_case_ab}
LU^1=-\quadre{a^2(x)H+\mu(x, t_1)}U^0_z.
\end{equation}
Multiplying by $U^0_z$ and integrating by parts, we obtain
\begin{equation*}
\int_{-\infty}^{+\infty}  LU^1 U^0_z dz = \int_{-\infty}^{+\infty}LU^0_z U^1 dz = 0.
\end{equation*}
Thus
\begin{equation*}
\int_{-\infty}^{+\infty}LU^1  U^0_z dz=-(a^2(x)H(p(x), t_1))+\mu(x,
t_1))\int_{-\infty}^{+\infty} (U^0_z)^2 dz=0.
\end{equation*}
Since $\int_{-\infty}^{+\infty} (U^0_z)^2 dz\neq 0$, it follows $\mu(x, t_1)=-a^2(x)H(p(x), t_1)$. Hence,
\begin{equation*}
W^1(z, x, t_1)=-a^2(x)H(p(x), t_1)U^0_z\left[\frac{b(x)}{a(x)}
z\right].
\end{equation*}
From (\ref{LU1_case_ab}), we have $LU^1=0$, so that $U^1(z, x,
t_1)=c(x, t_1)U^0_z\quadre{\frac{b(x)}{a(x)}z}$. 
Due
to (\ref{normalization_ab}) and to (\ref{match2}), the function $c(x, t_1)$ is null and therefore
\begin{equation*}
U^1(z, x, t_1)=0.
\end{equation*}
\subsubsection*{Second order}
At second order in $\eps$, using the definition of the operator $L$, we deal with the equation
\begin{eqnarray}\label{inn_two_case_ab1}
LW^2 & = & -a^2(x)W^1_z H-2 a(x) \nabla a(x) W^1_z \nabla d \nonumber\\  
         & = & - [a^2(x) H+2a(x) \nabla a(x)\nabla d]W^1_z. 
\end{eqnarray}
From (\ref{w_case_ab_app}), we have
\begin{eqnarray*}
W^2& = & b^2(x)\psi_{uu}(U^0)U^2-2 a(x) \nabla a(x)U^0_z
\nabla
d\\
& - & a^2(x) \Delta_x U^0+a^2(x) z U^0_z (H^2-2K)-a^2(x) U^2_{zz},
\end{eqnarray*}
which gives
\begin{eqnarray}\label{inn_two_case_ab2}
LU^2& = & -W^2-2 a(x) \nabla a(x) \nabla d U^0_z \nonumber\\
& - & a^2(x)\Delta_x U^0
+ a^2(x)z U^0_z (H^2-2K).
\end{eqnarray}
Calculations give the following expression for the laplacian of $ U^0$
\begin{equation*}
\Delta_x U^0=\frac{a(x)}{b(x)}\Delta\quadre{\frac{b(x)}{a(x)}} z U^0_z+\frac{a^2(x)}{b^2(x)}\quadre{ \nabla\tonde{ \frac{b(x)}{a(x)}}}^2z^2 U^0_{zz}.
\end{equation*}
Substituting it into \eqref{inn_two_case_ab2}, we obtain
\begin{eqnarray}\label{LU2}
LU^2 &=&  -W^2-2 a(x) \nabla a(x) \nabla d U^0_z\\
& + &  a^2(x)z U^0_z \graffe{H^2-2K -\frac{a(x)}{b(x)}\Delta\quadre{\frac{b(x)}{a(x)}}}\nonumber \\
& - & \frac{a^4(x)}{b^2(x)} \quadre{ \nabla\tonde{ \frac{b(x)}{a(x)}}}^2z^2 U^0_{zz}.\nonumber
\end{eqnarray}
In order to obtain the function $W^2$, we remind the matching condition (\ref{match3})
\begin{equation*}
\lim_{z\to \pm \infty} W^2(z, x, t_1)=0.
\end{equation*}
The solution of (\ref{inn_two_case_ab1}) is then
\begin{equation*}
W^2 = H^2(p(x), t_1)\quadre{a^2(x) \frac{z}{2}+\alpha(x, t_1)}U^0_z+a(x) H(p(x), t_1) \nabla a(x) \nabla d z U^0_z,
\end{equation*}
where $\alpha(x, t_1)$ is a bounded function.
Therefore
\begin{eqnarray*}
LU^2= & &  \graffe{ \frac{a^2(x) H^2}{2}-2 a^2(x) K-a(x) H \nabla a(x) \nabla d-\frac{a^3(x)}{b(x)}\Delta\quadre{\frac{b(x)}{a(x)}}}z U^0_z\\
& - &  \quadre{\alpha H^2 +2 a(x) \nabla a(x) \nabla d}U^0_z\\
& - &  a^2(x) \frac{a^2(x)}{b^2(x)}\quadre{\nabla\tonde{\frac{b(x)}{a(x)}}}^2 z^2 U^0_{zz}.
\end{eqnarray*}
We remind that $U^2$ has to satisfy the matching condition (\ref{match3})
\begin{equation*}
\lim_{z\to\pm \infty} U^2(z, x, t_1)=0.
\end{equation*}
With the same notation of \citeasnoun{LM}, if we set
\begin{equation}\label{U_2}
U^2(z, x, t_1)=f(z, x, t_1)U^0_z,
\end{equation}
we find that 
\begin{equation*}
f_z=\frac{g}{a^2(x) (U^0_z)^2},
\end{equation*}
where $g$ is such that
\begin{eqnarray*}
g_z & = &\graffe{ \frac{a^2(x) H^2}{2}-2 a^2(x) K-a(x) H \nabla a(x) \nabla d-\frac{a^3(x)}{b(x)}\Delta\quadre{\frac{b(x)}{a(x)}}}z (U^0_z)^2\\
& - & \quadre{\alpha H^2 +2 a(x) \nabla a(x) \nabla d} (U^0_z)^2\\
& + & 2\frac{a^3(x)}{b(x)}\quadre{\nabla \tonde{\frac{b(x)}{a(x)}}}^2 z^2U^0 (U^0_z)^2.
\end{eqnarray*}
\subsubsection*{Third order}
At third order in $\eps$, the first equation of problem \eqref{problem_ab} becomes
\begin{eqnarray}{\label{inn_third_case_ab}}
\frac12 VU^0_z  = & - & \frac12
W^1+LW^3-b^2(x)\psi_{uuu}(U^0)U^2W^1\\
& - & a^2(x)z W^1_z (H^2-2K)+a^2(x)W^2_z H+a^2(x) \Delta_xW^1.\nonumber
\end{eqnarray}
From now on, we omit to write the dependance on $x$ of $a$ and $b$.
Here we need the expansion of the term $\Delta_x W^1$.
It may be obtained calculating
\begin{equation*}
\Delta_xW^1 = \textrm{div}(\nabla W^1) =  -\textrm{div}\quadre{\frac{a}{b}\nabla(abH)U^0_z}+2 z
\textrm{div}\quadre{a^2\nabla \tonde{\frac{b}{a}}H U^0 U^0_z}.
\end{equation*}
Using \eqref{relazione} and collecting the different terms that we obtain from the calculations, we can write 
\begin{equation*}
\Delta_xW^1 = A z U^0_{zz}+B U_z^0+C z^2 (U_z^0)^2+D z^2 U^0 U^0_{zz},
\end{equation*}
where (from \cite{cen}, $\Delta_x H=\Delta_s H$)
\begin{eqnarray*}
A: & = &  2(\nabla a)^2 H - 2\frac{a^2}{b^2}(\nabla b)^2 H-2\frac{a^2}{b}\nabla b \nabla H+2a\nabla a \nabla H\\
& - & \frac{a^2}{b} \Delta b H+ a \Delta a H+2\frac{a}{b}\nabla a \nabla b H-2 \frac{(\nabla a)^2}{b}H,
\end{eqnarray*}
\begin{eqnarray*}
B:& = & -a H \Delta a-2a \nabla a \nabla H-2\frac{a}{b}
\nabla
a\nabla b H -\frac{a^2}{b}H \Delta b\\
& - & 2\frac{a^2}{b} \nabla b \nabla H-a^2 \Delta_s H,
\end{eqnarray*}
\begin{equation*}
C:=2 a \frac{(\nabla b)^2}{b} H-4 \nabla a \nabla b H-2 \frac{b}{a}(\nabla a)^2 H,
\end{equation*}
\begin{equation*}
D:=2\frac{a}{b} (\nabla b)^2 H-2\nabla a \nabla b H+2 \frac{b}{a}(\nabla a)^2 H.
\end{equation*}
Thus
\begin{eqnarray}\label{almost_V}
\frac12VU^0_z  = & & \frac12
a^2H U^0_z+LW^3+a^2 b^2H\psi_{uuu}(U^0)U^2U^0_z \nonumber\\
& + & a^4 H z U^0_{zz} (H^2-2K)+a^2W^2_z H \nonumber \\
& + & a^2 A z U^0_{zz}+a^2 B U^0_z+a^2 C z^2 (U^0_z)^2 + a^2 D z^2 U^0U^0_{zz}.
\end{eqnarray}
Taking into account the fact that $\psi_{uuu}(u)=12 u$, we multiply equation (\ref{almost_V}) by $U^0_z$ 
and integrate in $z$
\begin{eqnarray*}
\frac12V \int_{-\infty}^{+\infty} (U^0_z)^2 dz  & =  & \frac12
a^2H\int_{-\infty}^{+\infty}  (U^0_z)^2 dz+\int_{-\infty}^{+\infty} LW^3U^0_z dz\\
& + & 12a^2 b^2H \int_{-\infty}^{+\infty} U^0U^2(U^0_z)^2 dz\\
& + & a^4 H (H^2-2K) \int_{-\infty}^{+\infty} z U^0_z U^0_{zz} dz\\
& + & a^2 H \int_{-\infty}^{+\infty}W^2_z U^0_z dz+ a^2 A \int_{-\infty}^{+\infty} z U^0_z U^0_{zz}dz\\
&+ &a^2 B  \int_{-\infty}^{+\infty}(U^0_z)^2 dz+a^2 C  \int_{-\infty}^{+\infty}z^2 (U^0_z)^3 dz\\
& + & a^2 D  \int_{-\infty}^{+\infty} z^2 U^0U^0_zU^0_{zz}dz.
\end{eqnarray*}
Integration by parts yields
\begin{equation*}
\int_{-\infty}^{+\infty} LW^3U^0_z dz=0.
\end{equation*}
We set
\begin{equation*}
i_2=\int_{-\infty}^{+\infty} (U^0_z)^2 dz,
\end{equation*}
therefore
\begin{eqnarray}\label{V_int}
\frac12 V i_2 & = & \frac12a^2 Hi_2+12a^2 b^2 H \int_{-\infty}^{+\infty} U^0 U^2 (U^0_z)^2 dz\\
& + & a^4 H(H^2-2 K) \int_{-\infty}^{+\infty} z U^0_z U^0_{zz} dz\nonumber \\
& + & a^2 H \int_{-\infty}^{+\infty} U^0_z W^2_z dz\nonumber \\
& + & a^2 A \int_{-\infty}^{+\infty} z U^0_z U^0_{zz} dz+ a^2 B i_2\nonumber \\
& + & a^2 C \int_{-\infty}^{+\infty}z^2 (U^0_z)^3 dz\nonumber \\
& + & a^2 D \int_{-\infty}^{+\infty} z^2 U^0 U^0_z U^0_{zz}dz.
\end{eqnarray}
We calculate explicitly the integrals in the above expression.
From integration by parts, it is easy to obtain
\begin{equation*}
\int_{-\infty}^{+\infty} z U^0_z U^0_{zz}
dz = -\frac{i_2}{2}
\end{equation*}
and consequently
\begin{eqnarray*}
\int_{-\infty}^{+\infty} z U^0 (U^0_z)^2 dz =  -\frac12\frac{a}{b} \int_{-\infty}^{+\infty} z U^0_z U^0_{zz} dz=\frac14 \frac{a}{b}i_2.
\end{eqnarray*}
Afterwards, due to \eqref{relazione} and \eqref{U_2}
\begin{eqnarray*}
\int_{-\infty}^{+\infty} U^0 U^2 (U^0_z)^2
dz & = & - \frac12\frac{a}{b} \int_{-\infty}^{+\infty}U^0_{zz} f (U^0_z)^2 dz \\
& = & \frac16 \frac{1}{ab}\int_{-\infty}^{+\infty} g U^0_z dz\\
& = & - \frac16 \frac{1}{ab}\int_{-\infty}^{+\infty} g_z U^0 dz\\
& = & - \frac16 \frac{1}{ab}  \quadre{ \frac{a^2H^2}{2}-2 a^2 K-a H \nabla a \nabla d-\frac{a^3}{b}\Delta\tonde{\frac{b}{a}}}  \int_{-\infty}^{+\infty}z U^0(U^0_z)^2 dz\\
& + & \frac16 \frac{1}{ab}  \quadre{\alpha H^2 +2 a \nabla a \nabla d}\int_{-\infty}^{+\infty} U^0 (U^0_z)^2 dz\\
& - & \frac16 \frac{1}{ab}  2 \frac{a^3}{b}\quadre{\nabla \tonde{\frac{b}{a}}}^2\int_{-\infty}^{+\infty}  z^2(U^0)^2 (U^0_z)^2dz.
\end{eqnarray*}
The integral $\int_{-\infty}^{+\infty} U^0 (U^0_z)^2 dz$ is null since the integrand is an odd function with respect to $z$, while
\begin{equation*}
\int_{-\infty}^{+\infty} z^2 (U^0)^2 (U^0_z)^2 dz=\frac{\pi^2}{45}\frac{a^3}{b}. 
\end{equation*}
Thus
\begin{eqnarray*}
 \int_{-\infty}^{+\infty} U^0 U^2 (U^0_z)^2 dz & = &-\frac{1}{24}\frac{i_2}{b^2}  \quadre{ \frac{a^2H^2}{2}-2 a^2 K-a H \nabla a \nabla d-\frac{a^3}{b}\Delta\tonde{\frac{b}{a}}}\\
& - & \frac13 \frac{a^5}{b^3}\frac{\pi^2}{45}\quadre{\nabla\tonde{\frac{b}{a}}}^2.
\end{eqnarray*}
Furthermore
\begin{eqnarray*}
\int_{-\infty}^{+\infty} U^0_z W^2_z dz & = &\int_{-\infty}^{+\infty} U^0_z\quadre{\frac{H^2}{2}\tonde{a^2+\frac{2a}{H} \nabla a \nabla d}\tonde{U^0_z+ z U^0_{zz}}}dz\\
& + & \int_{-\infty}^{+\infty} \alpha H^2 U^0_z U^0_{zz} dz\\
& = & \frac{H^2}{2} \tonde{a^2+\frac{2a}{H} \nabla a \nabla d} \int_{-\infty}^{+\infty} (U^0_z)^2\\
& + & \frac{H^2}{2}\tonde{a^2+\frac{2a}{H} \nabla a \nabla d}\int_{-\infty}^{+\infty} z U^0_z U^0_{zz}\\& + &  \alpha H^2  \int_{-\infty}^{+\infty} U^0_z U^0_{zz} dz\\
& = & \tonde{\frac{a^2 H^2}{4}+\frac{aH \nabla a \nabla d}{2}}i_2.
\end{eqnarray*}
Finally
\begin{equation*}
\int_{-\infty}^{+\infty}z^2 (U^0_z)^3 dz=\frac{4 \pi^2}{45}-\frac23
\end{equation*}
and 
\begin{equation*}
\int_{-\infty}^{+\infty} z^2 U^0 U^0_z U^0_{zz}dz=-\frac{2
\pi^2}{45}.
\end{equation*}
The substitution of the above integrals into \eqref{V_final} gives the following interface equation
\begin{eqnarray*}
V & = & a^2 H+4 H a^4\tonde{K-\frac{H^2}{4}}\\
& + & 2 a^3 H^2 \nabla a \nabla d+\frac{a^5}{b}H\Delta\tonde{\frac{b}{a}}-\frac{2}{15}\pi^2 a^6 H \quadre{\nabla \tonde{\frac{b}{a}}}^2\\
& - & a^2 A+2a^2 B+abC\tonde{\frac{2}{15}\pi^2-1}-\frac{\pi^2}{15} a b D,
\end{eqnarray*}
if $N=3$.
\section*{Conclusions}
In this paper, we have studied a geometric law considering the gradient flow of the functional (\ref{func_ab}), built following the variational motivation used in \citeasnoun{LM}. We have proved, by means of formal asymptotics, that such a motion may be approximated in a suitable sense by an inhomogeneous fourth order parabolic equation. It involves the mean and Gaussian curvatures and the surface Laplacian of the mean curvature of the evolving interface. As in \citeasnoun{nmhs}, the motion equation arising from \eqref{u_t_case_ab_app} and \eqref{w_case_ab_app} involves drift terms, despite the absence of drifts in the original equation. 
\section*{Acknowledgements}
The author is grateful to Paola Loreti and Riccardo March, for  fruitful mathematical discussions about this subject.
\bibliographystyle{dcu}
\bibliography{References}

\end{document}